\documentclass[12pt]{amsart}
\usepackage{epsfig,amsfonts,amsthm,amssymb,latexsym,amsmath}
\usepackage{subfigure}
\usepackage{graphicx,color}

\theoremstyle{plain}

\theoremstyle{remark}

\theoremstyle{definition}

\title{CENTRAL SCHEMES FOR POROUS MEDIA FLOWS}

\vspace{5pt}

\author{\scriptsize E. \ Abreu,$^{1}$ \ 
\ F. \ Pereira$^{2}$ \ and \ S. \ Ribeiro$^{2}$}
\date{}

\begin{document}

\def \RR {{\mathbb R}}
\def \ZZ {{\mathbb Z}}
\def \NN {{\mathbb N}}
\def \vx {{^x\!v}}
\def \vy {{^y\!v}}
\def \dx {{\Delta X}}
\def \dy {{\Delta Y}}
\def \dt {{\Delta t}}
\def \xip {{x^n_{j+1/2,l}}}
\def \xfp {{x^n_{j+1/2,r}}}
\def \xim {{x^n_{j-1/2,l}}}
\def \xfm {{x^n_{j-1/2,r}}}
\def \yip {{y^n_{k+1/2,l}}}
\def \yfp {{y^n_{k+1/2,r}}}
\def \yim {{y^n_{k-1/2,l}}}
\def \yfm {{y^n_{k-1/2,r}}}
\def \tmk {{t^m_\kappa}}

\maketitle

\vspace{-20pt}
\begin{center}
{\footnotesize $^1$Instituto Nacional de Matem{\'a}tica Pura e Aplicada, RJ 22460-320, Brazil \\
$^2$Department of Mathematics and School of Energy Resources,\\ University of Wyoming, Laramie, WY 82071-3036, USA. \\
$^3$ Departamento de Ci\^encias Exatas, Universidade Federal Fluminense, Volta Redonda/RJ 27255-250, Brazil \\ 
E-mails: eabreu@impa.br / lpereira@uwyo.edu / sribeiro@metal.eeimvr.uff.br
}\end{center}

\hrule

\begin{abstract}
We are concerned with central differencing schemes for solving scalar hyperbolic
conservation laws arising in the simulation of multiphase flows in heterogeneous 
porous media. We compare the Kurganov-Tadmor (KT) \cite{KT2000} semi-discrete 
central scheme with the Nessyahu-Tadmor (NT) \cite{TN} central scheme. The KT scheme 
uses more precise information about the local speeds of propagation together with 
integration over nonuniform control volumes, which contain the Riemann fans. 
These methods can accurately resolve sharp fronts in the fluid saturations without 
introducing spurious oscillations or excessive numerical diffusion. We first discuss the coupling of
these methods with velocity fields approximated by mixed finite elements. Then,
numerical simulations are presented for two-phase, two-dimensional flow problems 
in multi-scale heterogeneous petroleum reservoirs. We find the KT scheme to be 
considerably less diffusive, particularly in the presence of high permeability flow channels, which 
lead to strong restrictions on the time step selection; however, the KT scheme may 
produce incorrect boundary behavior.
\end{abstract}

\medskip
\noindent

\medskip
\noindent
\medskip

\hrule

\section{Introduction}\label{sec1}

We are concerned with high resolution central schemes for solving 
scalar hyperbolic conservation laws arising in the simulation of multiphase 
flows in multidimensional heterogeneous petroleum reservoirs. 

Many of the modern high resolution approximations for
nonlinear conservation laws employ Godunov's appoach \cite{GOD1959}
or  \textbf{REA} (reconstruct, evolve, average)
algorithm, i.e., the approximate solution is represented by
a piecewise polynomial which is \textbf{R}econstructed from
the \textbf{E}volving cell \textbf{A}verages. The two main
classes of Godunov methods are upwind and central schemes.

The Lax-Friedrichs (LxF) scheme \cite{LAX1954} is the
canonical first order central scheme, which is the
forerunner of all central differencing schemes. It is based
on piecewise constant approximate solutions. It also enjoys
simplicity, i.e., it does not employ Riemann solvers and
characteristic decomposition. Unfortunately the excessive
numerical dissipation in the LxF recipe (of order 
${\mathcal O}\big((\Delta X)^{2}/\Delta t\big)$) yields poor
resolution, which seems to have delayed the development of
high resolution central schemes when compared with the
earlier developments of the high resolution upwind methods.
Only in 1990 a second order generalization to the LxF scheme
was introduced by Nessyahu and Tadmor (NT) \cite{TN}. They
used a staggered form of the LxF scheme and replaced the
first order piecewise constant solution with a van Leer's
MUSCL-type piecewise linear second order approximation
\cite{VANLEER1979}.  The numerical dissipation in this new
central scheme has an amplitude of order 
${\mathcal O}\big((\Delta X)^{4}/\Delta t\big)$. When applying these
methods to multiphase flows in highly heterogeneous
petroleum reservoirs or aquifers we need to use
decreasing time steps as the heterogeneity increases,
yielding  greater numerical diffusion. 
Kurganov and Tadmor (KT) \cite{KT2000} combined
ideas from the construction of the NT scheme with  Rusanov's
method \cite{RUSANOV} to obtain the first second order central scheme
that admits a semi-dis\-cre\-te formulation which is then solved with
an appropriate ODE solver. The resulting scheme has a much
smaller numerical diffusion than the NT scheme. Due to the semi-discrete
formulation, this numerical diffusion is independent of the time step
used to evolve the ordinary differential equation. This property guarantees 
that no extra numerical diffusion will be added if the time step
is forced to decrease. For this reason, the application of this central scheme results in a new numerical approach to model two-phase flows with a much lower numerical diffusion, even in the presence of a highly heterogeneous porous media.

The goals of this paper are (i) to discuss the coupling of NT and KT schemes to
velocity fields approximated by Raviart-Thomas mixed finite element method (See \cite{prta}), and
(ii) to compare the KT semi-discrete central scheme 
with the NT central scheme for numerical simulations of two-phase, incompressible,
two-dimensional flows in heterogeneous formations. Both methods can accurately 
resolve sharp fronts in the fluid saturations without introducing spurious 
oscillations or excessive numerical diffusion. 

Our numerical experiments indicate that the KT scheme is considerably less 
diffusive, particularly in the presence of viscous fingers, which lead to strong 
restrictions on the time step selection. On the other hand the KT scheme may 
produce incorrect boundary behavior in a typical two-dimensional geometry used 
in the study of porous media flows: the quarter of a five spot. 

Numerous methods have been introduced to solve two-phase
flow problems in porous media. Among eulerian-lagrangian
procedures we mention the Modified Method of Characteristics
\cite{PAPER_EWING,DouglasRussel1982}, the Modified Method of Characteristics
with Adjusted Advection \cite{dfp}, the Locally Conservative
Eulerian Lagrangian Method \cite {lcelm} and Eulerian
Lagragian Localized Adjoint Methods \cite{dahle}.
Additional techniques, to name just a few, include
higher--order Godunov schemes \cite{bell}, the
front-tracking method \cite{PAPER_GLIMM}, the streamline method
~\cite{streamlineSPE,streamline_thesis}
 the streamline
upwind Petrov-Galerkin method (SUPG)
\cite{loula_coutinho1994,coutinho1999} and a second-order
TVD-type finite volume scheme \cite{Durlofsky1993} 
(this procedure aims at the modeling of flow
through geometrically complex geological reservoirs).
Each of
these procedures has advantages and disadvantages.  We
refer the reader to \cite{RECENT_EWING,lcelm} and references
cited there for a discussion of these methods.

We remark that central schemes are particularly interesting for the numerical 
simulation of multiphase flow problems in porous media because they have been 
formulated to solve hyperbolic systems; this is not the case for several of the 
procedures mentioned above, which have been developed only for scalar equations. 

Moreover  these central schemes were also used to deal with many other applied problems: to solve Hamilton-Jacobi Equations (see \cite{Hamilton-Jacobi-NT} and \cite{Hamilton-Jacobi-KT}), to model the two-dimensional magnetohydrodynamics (MHD) equations and to study the Orszag-Tang vortex system, which describes the transition to supersonic turbulence for the equations of MHD in two space dimensions (see \cite{MHD-NT} and \cite{MHD-KT}), to mention just a couple of them.  

This paper is organized as follows. In Section \ref{sec2}, we discuss our strategy 
for solving numerically the model for two-phase, immiscible and incompressible 
displacement in heterogeneous porous media considered here. 
In Section \ref{sec3}, we discuss the application of central differencing 
schemes to porous media flows. In Section \ref{sec:Numerical_experiments}, we present the computational 
solutions for the model problem considered here and our conclusions.

\section{Numerical approximation of two-phase flows}\label{sec2}

We consider a model for two-phase immiscible and incompressible displacement in 
heterogeneous porous media. The governing equations are strongly nonlinear and lead to 
shock formation, and with or without diffusive terms they are of practical importance 
in petroleum engineering \cite{dwpb,CJ1986}. See also \cite{fpcross} and the references 
therein for recent studies for the scale-up problem for such equations. 

The conventional theoretical description of two-phase flow in 
a porous medium, in the limit of vanishing capillary pressure, is via Darcy's 
law coupled to the Buckley-Leverett equation. The two phases will be referred 
to as water and oil, and indicated by the subscripts $w$ and $o$, respectively.
Without sources or sinks and neglecting the effects of capillarity and gravity, these equations read (See \cite{dwpb} for more details)
\begin{equation}
\nabla \cdot {\bf v}=0, \quad {\bf v}=-\lambda(s) K({\bf x}) \nabla p,  
\label{preeq}
\end{equation} 
\begin{equation}
\frac{\partial s}{\partial t}+ \nabla \cdot (f(s){\bf v})=0,
\label{sateq}
\end{equation} 
Here, ${\bf v}$ is the total seepage velocity, $s$ is the water saturation, 
$K({\bf x})$ is the absolute permeability, and $p$ is the pressure. The constant 
porosity has been scaled out by a change of the time variable. 
The total mobility, $\lambda(s)$, and water fractional flow function, $f(s)$, are defined in terms of the relative permeabilities $k_{ri}(s)$ and phase viscosities $\mu_i$ by 
\[
\lambda(s)=\frac{k_{rw}(s)}{\mu_w}+\frac{k_{ro}(s)}{\mu_o}, \quad
f(s) = \frac{k_{rw}(s)/\mu_w}{\lambda(s)}.
\]

\subsection{Operator splitting for two-phase flow.}

An operator splitting technique is employed for the
computational solution of the saturation equation
(\ref{sateq}) and the pressure equation (\ref{preeq}) in
which they are solved sequentially with possibly distinct time
steps. This splitting scheme has proved to be
computationally efficient in producing accurate numerical
solutions for two-phase flows. We refer the reader to
\cite{dfp} and references therein for more details
on the operator splitting technique; see also~\cite{AFMPb,AFMP, ADFMP, abreu_splitting_2008} and \cite{Aquino_splitting_2008}   for applications of this strategy to three phase flows taking into account
capillary pressure (diffusive effects).

Typically, for computational efficiency larger time steps
are used for the pressure-velocity
calculation (Equation \ref{preeq}) than for the convection
calculation (Equation \ref{sateq}).  Thus, we introduce two time steps: $\Delta
t_c$ for the solution of the hyperbolic problem for
convection, and $\Delta t_p$ for the pressure-velocity
calculation so that $\Delta t_p \geq \Delta t_c$. We remark
that in practice variable time steps are always useful,
especially for the convection micro-steps subject
dynamically to a $CFL$ condition.

For the pressure solution  we use a (locally conservative) hybridized mixed
finite element discretization equivalent to cell-centered
finite differences \cite{fpcross,dfp}, which effectively
treats the rapidly changing permeabilities that arise from
stochastic geology and produces accurate velocity
fields. The pressure and Darcy velocity are
approximated at times $t^m = m\Delta t_p$,
$m=0,1,2,\dots$. The linear system resulting from the
discretized equations is solved by a preconditioned
conjugate gradient procedure (PCG) (See \cite{dfp} and the
references therein).  The saturation equation
is approximated at times $t^{m}_{\kappa}=t^{m} +
{\kappa}\Delta t_c$ for $t^{m} < t^{m}_{\kappa} \leq
t^{m+1}$.  We remark that we must specify the water
saturation at $t=0$.

\section{Central differencing schemes for porous media flows}\label{sec3}

In this section, we shall study the
family of high resolution, non-oscillatory, conservative
central differencing schemes introduced by Nessyahu and
Tadmor (NT)  and Kurganov and Tadmor (KT).
They will be applied  
to the numerical approximation of the
scalar hyperbolic conservation law modeling the convective
transport of fluid phases in two-phase flows. For the associated elliptic problem (Eq. \eqref{preeq}), we use the lowest order Raviart-Thomas \cite{prta}
locally conservative mixed finite elements.
These central schemes enjoy the main advantage of
Godunov-type central schemes: simplicity, i.e., they employ
neither characteristic decomposition nor approximate Riemann
solvers. This makes them universal methods that can be
applied to a wide variety of physical problems, including
hyperbolic systems. In the following sections
 we will discuss the
main ideas of the NT and KT central schemes coupled to the
mixed finite element discretization mentioned above. We will
not repeat here all the details involved in the development
of the NT and KT schemes; instead, we refer the reader to
\cite{TN} and \cite{KT2000} for this material.
 
\subsection{The Nessyahu-Tadmor scheme for two-phase flows}
\label{sec:NTscheme}

Consider the following scalar hyperbolic conservation law,
\begin{equation}
\frac{\partial s}{\partial t} + \frac{\partial }{\partial
  x}(\vx f(s)) +
\frac{\partial }{\partial y}(\vy f(s)) = 0,
\label{two2D}
\end{equation}
subject to prescribed initial data, $s(x,y,0) = S_0(x,y)$.
Here $\vx = \vx(x,y,t)$ and $\vy = \vy(x,y,t)$ denote the $x-$ and 
$y-$components of the velocity field ${\bf v}$ (see Eq. \ref{preeq}). 
To approximate \eqref{two2D} by the NT scheme, we begin with a piecewise constant solution
of the form $\sum_{j,k} \overline{S}^\kappa_{j,k}\chi_{j,k}(x,y)$, where $ \overline{S}^{\kappa}_{j,k} := \overline{S}(x_j,y_k,\tmk)$ is the approximate cell average at $t = \tmk$ associated with the cell $C_{j,k} = I_j \times I_k = [x_{j-1/2},x_{j+1/2}] \times [y_{k-1/2},y_{k+1/2}]$ and $\chi_{j,k}(x,y)$ is a characteristic function of the cell $C_{j,k}$.

We first reconstruct a piecewise linear approximation of the form
\begin{equation}
\begin{array}{lll}
s(x,y,\tmk) & = &\displaystyle{\sum_{j,k} \widetilde{S}^{\kappa}_{j,k}(x,y) \chi_{j,k}(x,y)} \\ \\
& = & \displaystyle \sum_{j,k}\left[\overline{S}^{\kappa}_{j,k}
+\frac{(S^{\kappa}_{j,k}\acute{)}}{\Delta X}(x-x_{j}) \right.
\\ \\& & \hspace{+1.0cm}
\left. \displaystyle+ \frac{(S^{\kappa}_{j,k}\grave{)}}{\Delta Y} (y-y_{k})\right] \chi_{j,k}(x,y) \\ \\ 
& & x_{j-1/2}\leq x\leq x_{j+1/2}, \quad 
y_{k-1/2}\leq y\leq y_{k+1/2}.
\end{array}
\label{inter_2D02}
\end{equation}
In Eq. (\ref{inter_2D02}), the discrete slopes along the $x$ and $y$ directions satisfy
\begin{subequations}
\label{eq:sxy}
\begin{eqnarray}
\frac{(S^{\kappa}_{j,k}\acute{)\,\,}}{\Delta X} =  
\frac{\partial}{\partial x}s(x=x_j, y=y_k,t^\kappa) + O(\Delta X) \\
\frac{(S^{\kappa}_{j,k}\grave{)\,\,}}{\Delta Y} = 
\frac{\partial}{\partial y}s(x=x_j, y=y_k,t^\kappa)+ O(\Delta Y),
\end{eqnarray}
\end{subequations}
 to guarantee second-order accuracy.

The reconstruction \eqref{inter_2D02} retains conservation, i.e.:
\begin{equation}
\hspace{-7 mm}\displaystyle\frac{1}{\Delta X\! \Delta Y} 
\displaystyle{\int_{x_{j-1/2}}^{x_{j+1/2}}}
\displaystyle{\int_{y_{k-1/2}}^{y_{k+1/2}}} \widetilde{S}^{\kappa}(x,y)\,dxdy = \overline{S}^{\kappa}_{j,k}.
\label{parte2D}
\end{equation}

\begin{figure}[h]
\centering
\begin{picture}(0,0)%
\includegraphics{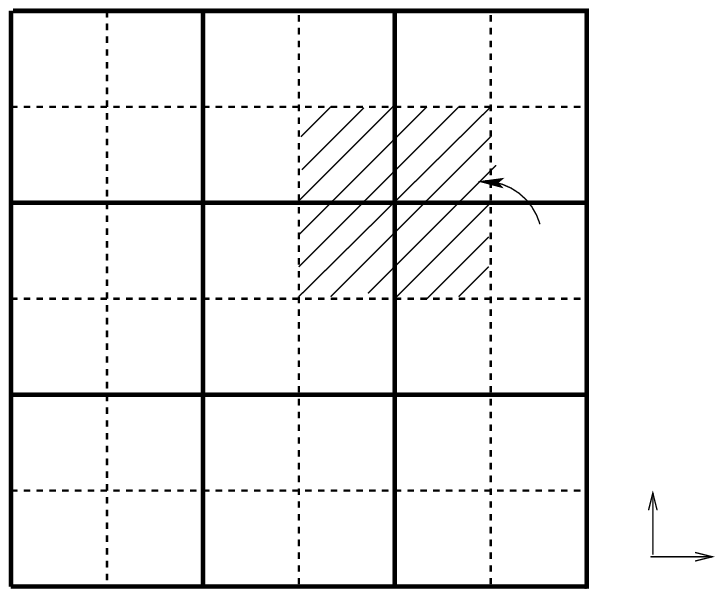}%
\end{picture}%
\setlength{\unitlength}{2693sp}%
\begingroup\makeatletter\ifx\SetFigFont\undefined%
\gdef\SetFigFont#1#2#3#4#5{%
  \reset@font\fontsize{#1}{#2pt}%
  \fontfamily{#3}\fontseries{#4}\fontshape{#5}%
  \selectfont}%
\fi\endgroup%
\begin{picture}(5232,4624)(3781,-5294)
\put(5851,-1546){\makebox(0,0)[lb]{\smash{\SetFigFont{8}{9.6}{\rmdefault}{\mddefault}{\updefault}{\color[rgb]{0,0,0}$(j,k+1)$}%
}}}
\put(4186,-1546){\makebox(0,0)[lb]{\smash{\SetFigFont{8}{9.6}{\rmdefault}{\mddefault}{\updefault}{\color[rgb]{0,0,0}$(j-1,k+1)$}%
}}}
\put(6931,-1546){\makebox(0,0)[lb]{\smash{\SetFigFont{8}{9.6}{\rmdefault}{\mddefault}{\updefault}{\color[rgb]{0,0,0}$(j+1,k+1)$}%
}}}
\put(3781,-5236){\makebox(0,0)[lb]{\smash{\SetFigFont{8}{9.6}{\rmdefault}{\mddefault}{\updefault}{\color[rgb]{0,0,0}$(x=0,y=0)$}%
}}}
\put(7861,-5206){\makebox(0,0)[lb]{\smash{\SetFigFont{8}{9.6}{\rmdefault}{\mddefault}{\updefault}{\color[rgb]{0,0,0}$(x=X,y=0)$}%
}}}
\put(8896,-4966){\makebox(0,0)[lb]{\smash{\SetFigFont{8}{9.6}{\rmdefault}{\mddefault}{\updefault}{\color[rgb]{0,0,0}$x$}%
}}}
\put(8371,-4381){\makebox(0,0)[lb]{\smash{\SetFigFont{8}{9.6}{\rmdefault}{\mddefault}{\updefault}{\color[rgb]{0,0,0}$y$}%
}}}
\put(3796,-826){\makebox(0,0)[lb]{\smash{\SetFigFont{8}{9.6}{\rmdefault}{\mddefault}{\updefault}{\color[rgb]{0,0,0}$(x=0,y=Y)$}%
}}}
\put(7831,-841){\makebox(0,0)[lb]{\smash{\SetFigFont{8}{9.6}{\rmdefault}{\mddefault}{\updefault}{\color[rgb]{0,0,0}$(x=X,y=Y)$}%
}}}
\put(4231,-3166){\makebox(0,0)[lb]{\smash{\SetFigFont{8}{9.6}{\rmdefault}{\mddefault}{\updefault}{\color[rgb]{0,0,0}$(j-1,k)$}%
}}}
\put(5851,-3181){\makebox(0,0)[lb]{\smash{\SetFigFont{8}{9.6}{\rmdefault}{\mddefault}{\updefault}{\color[rgb]{0,0,0}$(j,k)$}%
}}}
\put(6931,-3196){\makebox(0,0)[lb]{\smash{\SetFigFont{8}{9.6}{\rmdefault}{\mddefault}{\updefault}{\color[rgb]{0,0,0}$(j+1,k)$}%
}}}
\put(4246,-4546){\makebox(0,0)[lb]{\smash{\SetFigFont{8}{9.6}{\rmdefault}{\mddefault}{\updefault}{\color[rgb]{0,0,0}$(j-1,k-1)$}%
}}}
\put(5851,-4531){\makebox(0,0)[lb]{\smash{\SetFigFont{8}{9.6}{\rmdefault}{\mddefault}{\updefault}{\color[rgb]{0,0,0}$(j,k-1)$}%
}}}
\put(6976,-4531){\makebox(0,0)[lb]{\smash{\SetFigFont{8}{9.6}{\rmdefault}{\mddefault}{\updefault}{\color[rgb]{0,0,0}$(j+1,k-1)$}%
}}}
\put(7501,-2626){\makebox(0,0)[lb]{\smash{\SetFigFont{8}{9.6}{\rmdefault}{\mddefault}{\updefault}{\color[rgb]{0,0,0}$I_{j+\frac{1}{2},k+\frac{1}{2}}$}%
}}}
\end{picture}
\caption{Evolution step at each time level $t^{m}_{\kappa}$, $t^{m} < t^{m}_{\kappa} \leq t^{m+1}$, 
for the two-dimensional NT central differencing scheme.}
\label{grid_NT}
\end{figure}

Let $\{ s(x,y,t), t\ge \tmk\}$ be the exact solution of the conservation law \eqref{two2D}, subject to the reconstructed piecewise-linear data \eqref{inter_2D02} at time $t = \tmk$. The evolution step in the NT scheme consists of approximating this exact solution at the next time step $t = \tmk + \Delta t_c$, by its averages over staggered cells, $C_{j+1/2,k+1/2}  := I_{j+1/2} \times I_{k+1/2}$. See 
dashed grid in Figure \ref{grid_NT} (denote $\kappa + 1 := \tmk + \Delta t_c$).  Let
\begin{eqnarray}
\label{eq:cell_staggered}
\overline{S}^{\kappa +1}_{j+\frac{1}{2},k+\frac{1}{2}} & = & \displaystyle\frac{1}{\Delta X \! \Delta Y} \int_{C_{j+1/2,k+1/2}}  \hspace{-1.0cm}
s(x,y,t^{m}_{\kappa}+\Delta t_c)\,dxdy \nonumber
\\
 & = &
\displaystyle\frac{1}{\Delta X \! \Delta Y} \int_{x_j}^{x_{j+1}}\int_{y_k}^{y_{k+1}} 
s(x,y,t^{m}_{\kappa}+\Delta t_c)\,dxdy. \\
 & & x_{j} \leq x \leq x_{j+1}, \,\, y_{k} \leq y \leq y_{k+1} \nonumber
\end{eqnarray}
These new staggered cell averages are obtained by integrating the
conservation law (\ref{two2D}) over the control volumes $C_{j+1/2,k+1/2} \times [t^{m}_{\kappa},t^{m}_{\kappa}+\Delta t_c]$ following the same
manipulations as described in~\cite{NT2D} (denote 
$\alpha_x \equiv \frac{\Delta t_c}{\Delta X}$ and $\alpha_y \equiv \frac{\Delta t_c}{\Delta Y}$):
\begin{align}
\label{eq:stagerred_solution_future}
&\overline{S}^{\kappa +1}_{j+1/2,k+1/2}  =   
 \frac{1}{\Delta X \! \Delta Y} \int_{C_{j+1/2,k+1/2}}  \hspace{-1.0cm}
s(x,y,t^{m}_{\kappa}+\Delta t_c)\,dxdy \nonumber \\
 &-  \frac{\alpha_x}{\Delta X \! \Delta Y} \Big\{ \int_{\tmk}^{\tmk + \Delta t_c} \hspace{-0.3cm}\int_{y_k}^{y_{k+1}}  \hspace{-0.4cm}\Big[ \vx(x_{j+1},y,\tau)\, f(s(x_{j+1},y,\tau)  \nonumber \\
 &  \hspace{+4.5cm} - \vx(x_{j},y,\tau)\,f(x_j,y,\tau)\Big] dy\,d\tau \Big\} \nonumber \\ 
 &-  \frac{\alpha_y}{\Delta X \! \Delta Y} \Big\{ \int_{\tmk}^{\tmk + \Delta t_c} \hspace{-0.3cm}\int_{x_j}^{x_{j+1}}  \hspace{-0.4cm}\Big[\vy(x,y_{k+1},\tau)\, f(s(x,y_{k+1},\tau) \nonumber \\
 & \hspace{+4.5cm}- \vy(x,y_k,\tau)\, f(x,y_k,\tau)\Big] dx\,d\tau \Big\}.
\end{align}

The cell average $\int_{_{C_{j+1/2,k+1/2}}}  \hspace{-1.0cm} s(x,y,t^{m}_{\kappa}+\Delta t_c)\,dxdy$ has contributions from the four cells $C_{j,k}$, $C_{j+1,k}$, $C_{j+1,k+1}$, and $C_{j,k+1}$:
\begin{eqnarray}
\int_{C_{j+1/2,k+1/2}}  \hspace{-1.6cm}
s(x,y,t^{m}_{\kappa})\,dxdy &= & \int_{C_{j+1/2,k+1/2} \cap C_{j,k}}  \hspace{-2.0cm}
 \widetilde{S}^{\kappa}_{j,k}(x,y) +  \int_{C_{j+1/2,k+1/2} \cap C_{j,k+1}}  \hspace{-2.0cm}
 \widetilde{S}^{\kappa}_{j,k+1}(x,y) \nonumber \\
 & & +  \int_{C_{j+1/2,k+1/2} \cap C_{j+1,k}}  \hspace{-2.3cm}
 \widetilde{S}^{\kappa}_{j+1,k}(x,y) +  \int_{C_{j+1/2,k+1/2} \cap C_{j+1,k+1}}  \hspace{-2.5cm}
 \widetilde{S}^{\kappa}_{j+1,k+1}(x,y) 
\end{eqnarray}
Computing these integrals exactly yields 
\begin{eqnarray}
\label{eq:stagerred_solution_tn}
\overline{S}^{\kappa}_{j+\frac{1}{2},k+\frac{1}{2}} & = &
{\frac{1}{4}}(
   \overline{S}^{\kappa}_{j,k}
 + \overline{S}^{\kappa}_{j,k+1}
 + \overline{S}^{\kappa}_{j+1,k} + 
   \overline{S}^{\kappa}_{j+1,k+1}) \nonumber \\
   & & +  {\frac{1}{16}} \Big[
  (S^{\kappa}_{j,k}\acute{)\,\,} +  (S^{\kappa}_{j,k+1}\acute{)\,\,}
  -  (S^{\kappa}_{j+1,k}\acute{)\,\,} -  (S^{\kappa}_{j+1,k+1}\acute{)\,\,} \nonumber \\
& & \hspace{+1.0cm}+  
(S^{\kappa}_{j,k}\grave{)\,\,} -  (S^{\kappa}_{j,k+1}\grave{)\,\,}
  +  (S^{\kappa}_{j+1,k}\grave{)\,\,} -  (S^{\kappa}_{j+1,k+1}\grave{)\,\,}\Big].
\end{eqnarray}

To approximate the four flux integrals on the right hand side  of \eqref{eq:stagerred_solution_future}, we  use the second-order rectangular quadrature rule for the spatial integration and the midpoint quadrature rule for second-order approximation of the temporal integrals. For instance, letting  $\kappa + 1/2$ be $\tmk + \Delta t_c/2$,
\begin{subequations}
\label{eq:fluxes_integrals}
\begin{align}
 \frac{\alpha_x}{\Delta X \! \Delta Y}&  \int_{\tmk}^{\tmk + \Delta t_c} \hspace{-0.3cm} \int_{y_k}^{y_{k+1}}  \hspace{-0.4cm} \vx(x_{j+1},y,\tau) \,f(s(x_{j+1},y,\tau)) dy\,d\tau \approx \nonumber \\
&\approx \frac{\alpha_x}{2}\left[ \vx^{\kappa +1/2}_{j+1,k} \,f(s^{\kappa +1/2}_{j+1,k}) + \vx^{\kappa +1/2}_{j+1,k+1}\,f(s^{\kappa +1/2}_{j+1,k+1})\right],
\end{align}
 \begin{align}
 \frac{\alpha_y}{\Delta X \! \Delta Y} &  \int_{\tmk}^{\tmk + \Delta t_c} \hspace{-0.3cm}\int_{x_j}^{x_{j+1}}  \hspace{-0.4cm} \vy(x,y_{k+1},\tau) \, f(s(x,y_{k+1},\tau)) dy\,d\tau \approx \nonumber \\
 & \approx \frac{\alpha_y}{2}\left[\vy^{\kappa +1/2}_{j,k+1}\,f(s^{\kappa +1/2}_{j,k+1}) + \vy^{\kappa +1/2}_{j+1,k+1}\,f(s^{\kappa +1/2}_{j+1,k+1})\right].
\end{align}
\end{subequations}
Since these midvalues are computed at the center of  the cells, $C_{j,k}$, 
where the solution is smooth, provided an appropriate CFL condition is observed, we can use 
Taylor expansion together with the conservation law \eqref{two2D} to get
\begin{equation}
\label{eq:midvalues}
s^{\kappa +1/2}_{j,k} = \overline{S}^\kappa_{j,k}
-\frac{\alpha_x}{2}
 {\vx}^{\kappa}_{j,k} (f^{\kappa}_{j,k}\acute{)\,\,} 
-\frac{\alpha_y}{2} {\vy}^{\kappa}_{j,k}  (f^\kappa_{j,k}\grave{)\,\,}.
\end{equation}
Here, $(f^{\kappa}_{j,k}\acute{)\,\,}$ and $(f^{\kappa}_{j,k}\grave{)\,\,}$ are one-dimensional discrete slopes in the $x$ and $y$ directions, respectively. They satisfy the conditions
\begin{subequations}
\label{eq:fxy}
\begin{eqnarray}
\frac{(f^{\kappa}_{j,k}\acute{)\,\,}}{\Delta X} =  
\frac{\partial}{\partial x}f(s(x=x_j, y=y_k,t^\kappa)) + O(\Delta X) \\
\frac{(f^{\kappa}_{j,k}\grave{)\,\,}}{\Delta Y} = 
\frac{\partial}{\partial y}f(s(x=x_j, y=y_k,t^\kappa))+ O(\Delta Y),
\end{eqnarray}
\end{subequations}
in order  to produce a
second order scheme for the approximation of (\ref{two2D}).
To avoid spurious oscillations, it is essential to reconstruct
the discrete derivatives given by Equations \eqref{eq:sxy} and \eqref{eq:fxy} with built-in nonlinear
limiters. In this work we use the following MinMod limiter
\begin{subequations}
 \label{eq:minmod_2D}
\begin{eqnarray}
  (S_x)^\kappa_{j,k} & \approx & \text{MM}\theta\frac{1}{\Delta x}\left\{\overline{S}^\kappa_{j-1,k},\overline{S}^\kappa_{j,k},
    \overline{S}^\kappa_{j+1,k}\right\} \nonumber \\
  & := & \text{MM}\left(\theta \frac{\Delta S^\kappa_{j+1/2,k}}{\Delta x}, 
  \frac{\Delta S^\kappa_{j-1/2,k} -\Delta S^\kappa_{j+1/2,k}}{2\Delta
    x}, \theta \frac{\Delta S^\kappa_{j-1/2,k}}{\Delta x}\right); \\
(f_x)^\kappa_{j,k} & \approx & \text{MM}\theta\frac{1}{\Delta x}\left\{{f}^\kappa_{j-1,k},{f}^\kappa_{j,k},
    {f}^\kappa_{j+1,k}\right\} \nonumber \\
  & := & \text{MM}\left(\theta \frac{\Delta f^\kappa_{j+1/2,k}}{\Delta x}, 
  \frac{\Delta f^\kappa_{j-1/2,k} -\Delta f^\kappa_{j+1/2,k}}{2\Delta
    x}, \theta \frac{\Delta f^\kappa_{j-1/2,k}}{\Delta y}\right),
\end{eqnarray}
\end{subequations}
where $\Delta$ is the centered difference, $\Delta
S^\kappa_{j+1/2,k} = \overline{S}^\kappa_{j+1,k} - \overline{S}^\kappa_{j,k}$.
We
refer the reader to \cite{TN} and \cite{KT2000} and the references therein
for the various options for the form of such discrete
derivatives.

In our sequential scheme, when solving for the saturation in time, the total velocity 
$\bf v$ is given by the solution of the velocity-pressure equation. 
Recall that the solution of Eq. \eqref{preeq} is approximated the lowest order Raviart-Thomas 
mixed finite element method. Thus, the computed total velocity ${\bf v}$ is 
discontinuous at the vertices of the original non-staggered grid cells. This constitutes 
a difficulty for the staggered scheme \eqref{eq:stagerred_solution_future}, which requires the values of the 
total velocity ${\bf v}$ at these vertices at every other time step. To avoid this 
difficulty we use the non-staggered version of the NT scheme.

To turn the staggered scheme \eqref{eq:stagerred_solution_future} into a non-staggered scheme, we re-average the
reconstructed values of the underlying staggered scheme, thus recovering the cell 
averages of the central scheme over the original non-staggered grid cells. First we 
reconstruct a piecewise bilinear interpolant at the time step $\kappa +1 := t^{m}_{\kappa}+\Delta t_c$
\begin{equation}
\begin{array}{lll}
\widetilde{S}^{\kappa+1}_{j+1/2,k+1/2}(x,y) 
& = & \displaystyle{\overline{S}^{\kappa+1}_{j+1/2,k+1/2}
+\frac{(S^{\kappa + 1}_{j+1/2,k+1/2}\acute{)}}{\Delta X}(x-x_{j+1/2})}
\\ & &
 \displaystyle{+ \frac{(S^{\kappa + 1}_{j+1/2,k+1/2}\grave{)}}{\Delta Y} (y-y_{k+1/2})} \\ \\
& & x_{j}\leq x\leq x_{j+1}, \quad 
y_{k}\leq y\leq y_{k+1},
\end{array}
\label{eq:reconstruct_future}
\end{equation}
as in Equation \eqref{inter_2D02}, through the staggered cell averages given by 
(\ref{eq:stagerred_solution_future}), and re-average it over the original grid cells, giving the 
following non-staggered scheme:
\begin{equation}
\begin{array}{ll}
\overline{S}^{\kappa + 1}_{j,k}
& = \displaystyle \frac{1}{4}\,
(\overline{S}^{\kappa + 1}_{j-1/2,k-1/2}
+\overline{S}^{\kappa + 1}_{j-1/2,k+1/2}
+\overline{S}^{\kappa + 1}_{j+1/2,k-1/2}
+\overline{S}^{\kappa + 1}_{j+1/2,k+1/2}) \\ 
& +\displaystyle \frac{1}{16}\,
\Big[(S^{\kappa+1}_{j-1/2,k-1/2}\acute{\,\,)} + (S^{\kappa+1}_{j-1/2,k+1/2}\acute{)\,\,} \nonumber \\
& \hspace{+1.5cm}- (S^{\kappa+1}_{j+1/2,k-1/2}\acute{)\,\,} - (S^{\kappa+1}_{j+1/2,k+1/2}\acute{)\,\,}\,\Big] \\ 
& +\displaystyle{\frac{1}{16}}\,
\Big[(S^{\kappa+1}_{j-1/2,k-1/2}\grave{)\,\,} - (S^{\kappa+1}_{j-1/2,k+1/2}\grave{)\,\,} \nonumber \\ 
&\hspace{+1.5cm} + (S^{\kappa+1}_{j+1/2,k-1/2}\grave{)\,\,} - (S^{\kappa+1}_{j+1/2,k+1/2}\grave{)\,\,}\,\Big].
\end{array}
\label{JT4.7}
\end{equation}

\subsection{The Kurganov-Tadmor  scheme for two-phase 
flows}
\label{subsec4.2}

The first multidimensional extension of the KT scheme 
was presented in  \cite{KT2000}. This extension used the 
dimension by dimension approach, that is, the numerical fluxes
computed along the $x$ and $y$ directions are viewed as a generalization
of the one-spatial-dimension numerical fluxes. This approach consists of the following steps: at each time step $\tmk$ and at each cell $I_{j,k}$,
\begin{enumerate}
\item [(i)] Compute the difference of the one-dimensional numerical flux in one spatial dimension  in the 
  $x$ direction keeping $y$  constant and equal to $y_k$.  Denote this difference by
  \begin{equation*}
    \label{eq:fluxo_x}
    \mathcal{F}^x_{j+1/2,k}(t) := \frac{H^x_{j+1/2,k}(t) - H^x_{j-1/2,k}(t)}{\Delta
    X}.
  \end{equation*}
  The numerical flux $H^x_{j+1/2,k}(t)$ is
  \begin{eqnarray}
H^x_{j+1/2,k}(t) & := & 
  \frac{1}{2}\Big[\vx_{j+1/2,k}(t)\, f(S^+_{j+1/2,k}(t)) \nonumber \\
  & & \hspace{+1.0cm} + 
  \vx_{j+1/2,k}(t)\, f(S^-_{j+1/2,k}(t))\Big]  \nonumber \\ 
& & - \frac{a^x_{j+1/2,k}(t)}{2}\left[S^+_{j+1/2,k}(t) -
  S^-_{j+1/2,k}(t)\right],  \label{eq:flux_dir_x}
  \end{eqnarray}
  where 
  \begin{eqnarray}
  S^+_{j+1/2,k}(t) & = & \widetilde{S}_{j+1,k}(x_{j+1/2},\mathbf{y_k},t) \nonumber \\
  & &  =\overline{S}_{j+1,k}(t) -  \frac{\Delta X}{2}(S_x)_{j+1,k}(t) \quad\mbox{and} \nonumber  \\
  S^-_{j+1/2,k}(t) & = & \widetilde{S}_{j,k}(x_{j+1/2},\mathbf{y_k},t) \nonumber \\
  & & = \overline{S}_{j,k}(t) + \frac{\Delta X}{2}(S_x)_{j,k}(t) 
  \end{eqnarray}
  are the corresponding right and left intermediate values of $\widetilde{S}(x,t^\kappa)$ at $(x_{j+1/2},y_k)$.

The local speed of wave propagation $a^x_{j+1/2,k}(t)$ is estimated at the cell boundaries $(x_{j+1/2},y_k)$ as the upper bound  
\begin{eqnarray}
  \label{eq:velo_prop_x}
  a^x_{j+1/2,k}(t) & \hspace{-0.3cm}= &\hspace{-0.3cm} \max_{\omega} \left\{|\vx_{j+1/2,k}(t)\,
  f'(\omega)|\right\},
\end{eqnarray}
where $\omega$ is a value between $S^+_{j+1/2,k}(t)$ and $S^-_{j+1/2,k}(t)$.
The velocity field used in the KT scheme is obtained
directly from the Raviart-Thomas space on the cell edges:
\begin{eqnarray*}
\vx_{j+1/2,k}(t) := {{(v_r)}}_{jk}(t), &  &  
\vx_{j-1/2,k}(t) := {{(v_l)}}_{jk}(t), 
\end{eqnarray*}
where $v_r$ and $v_l$ stand for the velocity on the
one-dimensional ``right'' and ``left'' faces of the cells.

\item [(ii)]Analogously,   compute the difference of the one-dimensional numerical flux in the $y$ direction keeping $x$  constant and equal to $x_j$. This difference is denoted by
  \begin{equation*}
    \label{eq:fluxo_y}
    \mathcal{F}^y_{j,k+1/2}(t) := \frac{H^y_{j,k+1/2}(t) - H^y_{j,k+1/2}(t)}{\Delta
    Y}.
  \end{equation*}
  The one dimensional numerical flux in the $y$ direction is
 \begin{eqnarray}
H^y_{j,k+1/2}(t) & := &
\frac{1}{2}\Big[\vy_{j,k+1/2}(t)\, f(S^+_{j,k+1/2}(t)) \nonumber \\
& & \hspace{+1.0cm} +
  \vy_{j,k+1/2}(t)\, f(S^-_{j,k+1/2}(t))\Big]  \nonumber \\ 
& & - \frac{a^y_{j,k+1/2}(t)}{2}\left[S^+_{j,k+1/2}(t) -
  S^-_{j,k+1/2}(t)\right].\label{eq:flux_dir_y} 
\end{eqnarray}
In a similar way, the correspondent ``up'' and ``down'' intermediate values of $\widetilde{S}(x,t^\kappa)$ at $(x_{j},y_{k+1/2})$ are 
\begin{eqnarray*}
S^+_{j,k+1/2}(t) & = & \overline{S}_{j,k+1}(t) - \frac{\Delta Y}{2}(S_y)_{j,k+1}(t) \quad\mbox{and} \\
 S^-_{j,k+1/2}(t) & = &  \overline{S}_{j,k}(t) + \frac{\Delta Y}{2}(S_y)_{j,k}(t).
 \end{eqnarray*}
The local speed of wave propagation $a^y_{j,k+1/2}(t)$ in the $y$ direction is estimated at the cell boundaries $(x_{j},y_{k+1/2})$ as the upper bound  
\begin{eqnarray}
  \label{eq:velo_prop_x}
  a^y_{j+1/2,k}(t) & \hspace{-0.3cm}= &\hspace{-0.3cm} \max_{\omega} \left\{|\vy_{j+1/2,k}(t)\,
  f'(\omega)|\right\},
\end{eqnarray}
where $\omega$ is a value between $S^+_{j,k+1/2}(t)$ and $S^-_{j,k+1/2}(t)$.
Analogously the velocity field in the $y$ direction is obtained 
directly from the Raviart-Thomas space on the cell edges:
\begin{eqnarray*}
\vx_{j,k+1/2}(t) := {{(v_u)}}_{jk}(t), &  &  
\vx_{j,k-1/2}(t) := {{(v_d)}}_{jk}(t), 
\end{eqnarray*}
where $v_u$ and $v_d$ stand for the velocity on the
``upper'' and ``lower'' faces of the cells.

\item [(iii)] The cell average $\overline{S}^{\kappa + 1}_{j,k}$ in the next
time step $\tmk + \Delta t_c$ is then the solution of the following differential equation
  \begin{eqnarray}
    \label{eq:sol_futuro_dxd}
    \frac{d}{dt}S_{j,k}(t) & = & - (F^x_{j+1/2,k}(t) + F^y_{j,k+1/2}(t)) \nonumber \\
    & = & - \frac{H^x_{j+1/2,k}(t) - H^x_{j-1/2,k}(t)}{\Delta
    X} - \frac{H^y_{j,k+1/2}(t) - H^y_{j,k+1/2}(t)}{\Delta Y},
  \end{eqnarray}
 \end{enumerate}

The numerical derivatives are computed using the MinMod limiter given by Equation
\eqref{eq:minmod_2D}. In our numerical experiments, the parameter $\theta$ assumes values $1 < \theta < 1.8$.

The two-dimensional semi-discrete formulation
\eqref{eq:sol_futuro_dxd} comprises a system of nonlinear ordinary
differential equations for the discrete unknows $\{S_{j,k}(t)\}$. To solve it, we 
integrate in time  introducing a variable time step $\Delta t^n$. 
Although the forward Euler scheme can be used, 
it may be advantageous to use higher order discretizations in numerical simulations. The numerical examples
presented below use third-order Runge-Kutta ODE solvers based on convex combinations
of forward Euler steps. See ~\cite{Shu1988} and \cite{ShuOsher} for more details on a whole family of such schemes.

\section{TWO-DIMENSIONAL NUMERICAL EXPERIMENTS}
\label{sec:Numerical_experiments}

We present and compare the results of numerical simulations
of two-dimensional, two-phase flows associated with two
distinct flooding problems using the KT and NT schemes.

In all simulations, the reservoir contains initially $79\%$
of oil and $21\%$ of water.  Water is injected at a constant
rate of $0.2$ pore volumes every year. The viscosity of oil
and water used are $\mu_o=10.0\,cP$ and
$\mu_w=0.05\,cP$.
 The relative permeabilities are assumed to be: 
\[
k_{ro}(s) = (1 - (1-s_{ro})^{-1}s)^{2}, \quad
k_{rw}(s) = (1-s_{rw})^{-2}(s-s_{rw})^{2},
\]
where $s_{ro}=0.15$ and $s_{rw}=0.2$ are the residual oil
and water saturations, respectively.

For the heterogeneous reservoir studies we consider a scalar
absolute permeability field $K({\bf x})$ taken to be
log-normal (a fractal field, see \cite{GLPZ} and
\cite{fpcross} for more details) with moderately large
heterogeneity strength. The spatially variable permeability
field is defined on a $256 \times 64$ grid with three
different values of the coefficient of variation $CV$ ($CV = 0.5$, $1.2$, $2.2$)
given by the ratio between the standard deviation and the mean value of the
permeability field.

We now discuss the simulations in the slab geometry. 
We consider two-dimensional flows in a rectangular,
heterogeneous reservoir ({\it slab geometry}) having $256$m
$\times$ $64$m with three different sizes of computational grid:
$256\times 64$, $512 \times 128$ and $1024 \times 256$ cells.
The boundary conditions and injection and production
specifications for the two-phase flow equations
(\ref{preeq})-(\ref{sateq}) are as follows. The injection is
made uniformly along the left edge of the reservoir and
the production is taken along the right edge; no flow is allowed 
along the edges appearing at the top and bottom of the reservoir.

\begin{figure}[htbp!]
\includegraphics[scale=0.45]{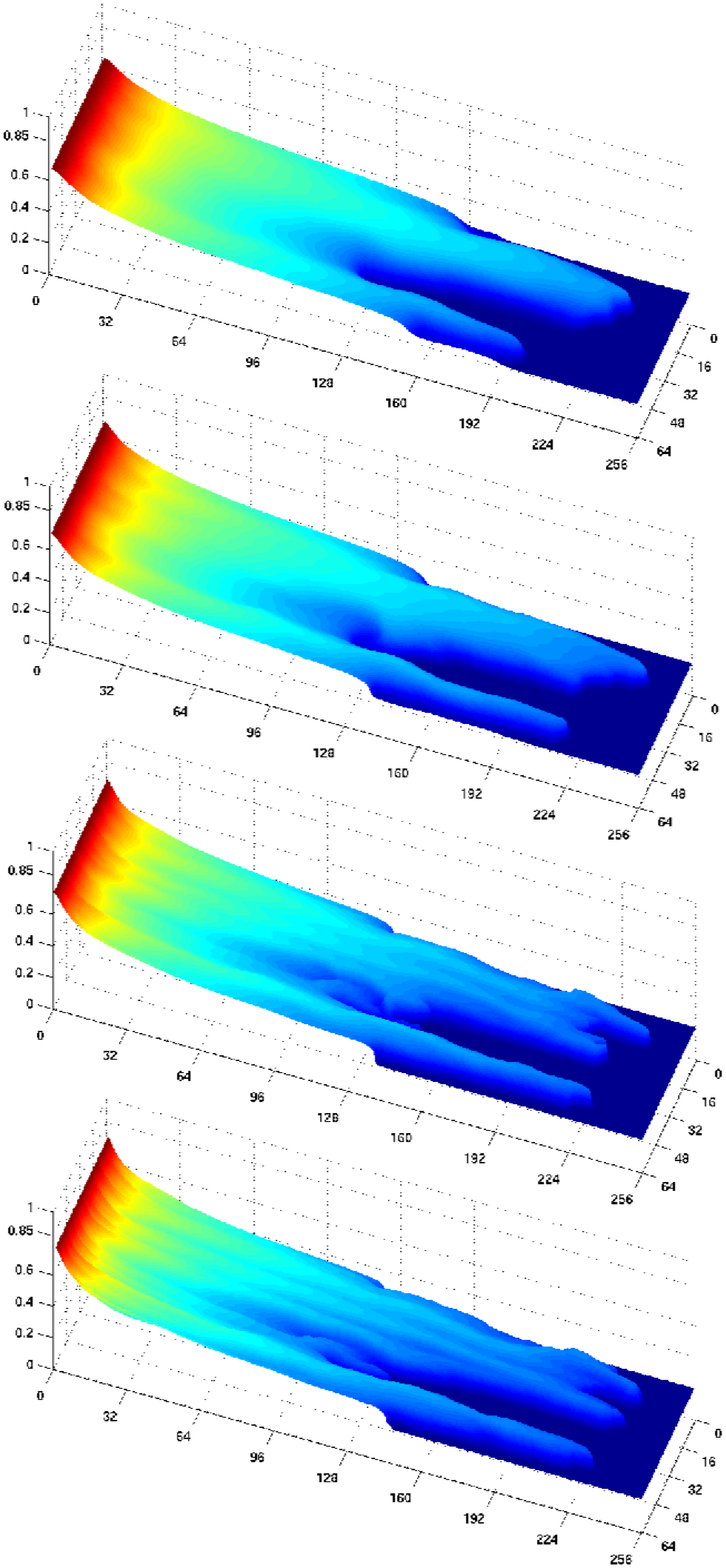}
\caption{Water saturation surface plots
  after $275$ days of simulation in a
  heterogeneous reservoir having $256$ m $\times$ $64$ m,
  with $CV=0.5$ and viscosity
  ratio 20. The first three pictures from top to bottom used
  the NT scheme with grids
  having $256 \times 64$, $512 \times 128$ and $1024 \times
  256$ cells, repectively. The bottom picture shows the KT scheme with a
  grid of $256 \times 64$ cells.}
  \label{slab}
\end{figure}

\begin{figure}[htbp!]
\includegraphics[scale=0.45]{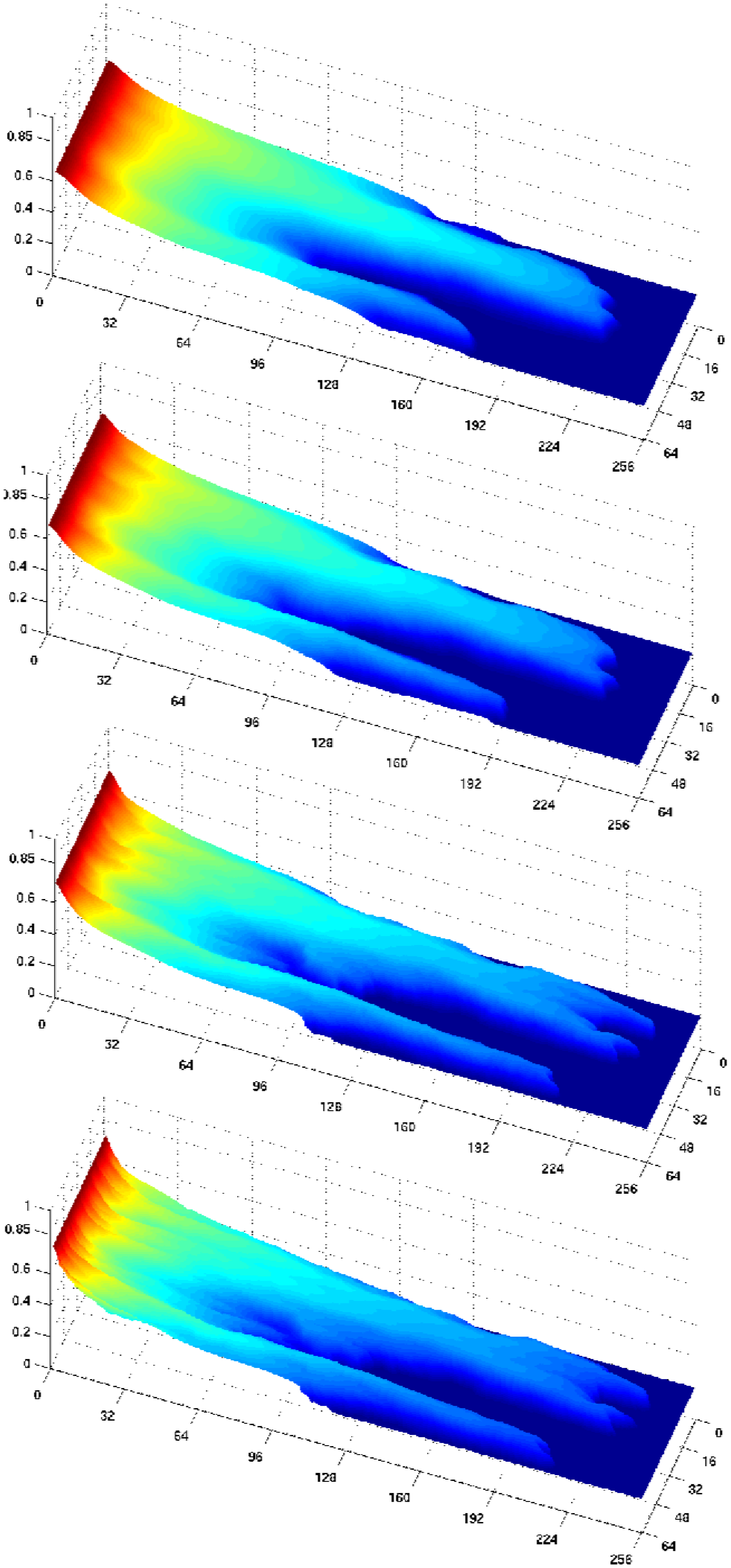}
\caption{Water saturation surface plots
  after $250$ days of simulation in a
  heterogeneous reservoir having $256$ m $\times$ $64$ m,
  with $CV=1.2$ and viscosity
  ratio 20. The first three pictures from top to bottom used
  the NT scheme with grids
  having $256 \times 64$, $512 \times 128$ and $1024 \times
  256$ cells, repectively. The bottom picture shows the KT scheme with a
  grid of $256 \times 64$ cells.}
\label{slab_CV1}
\end{figure}

\begin{figure}[htbp!]
\includegraphics[scale=0.45]{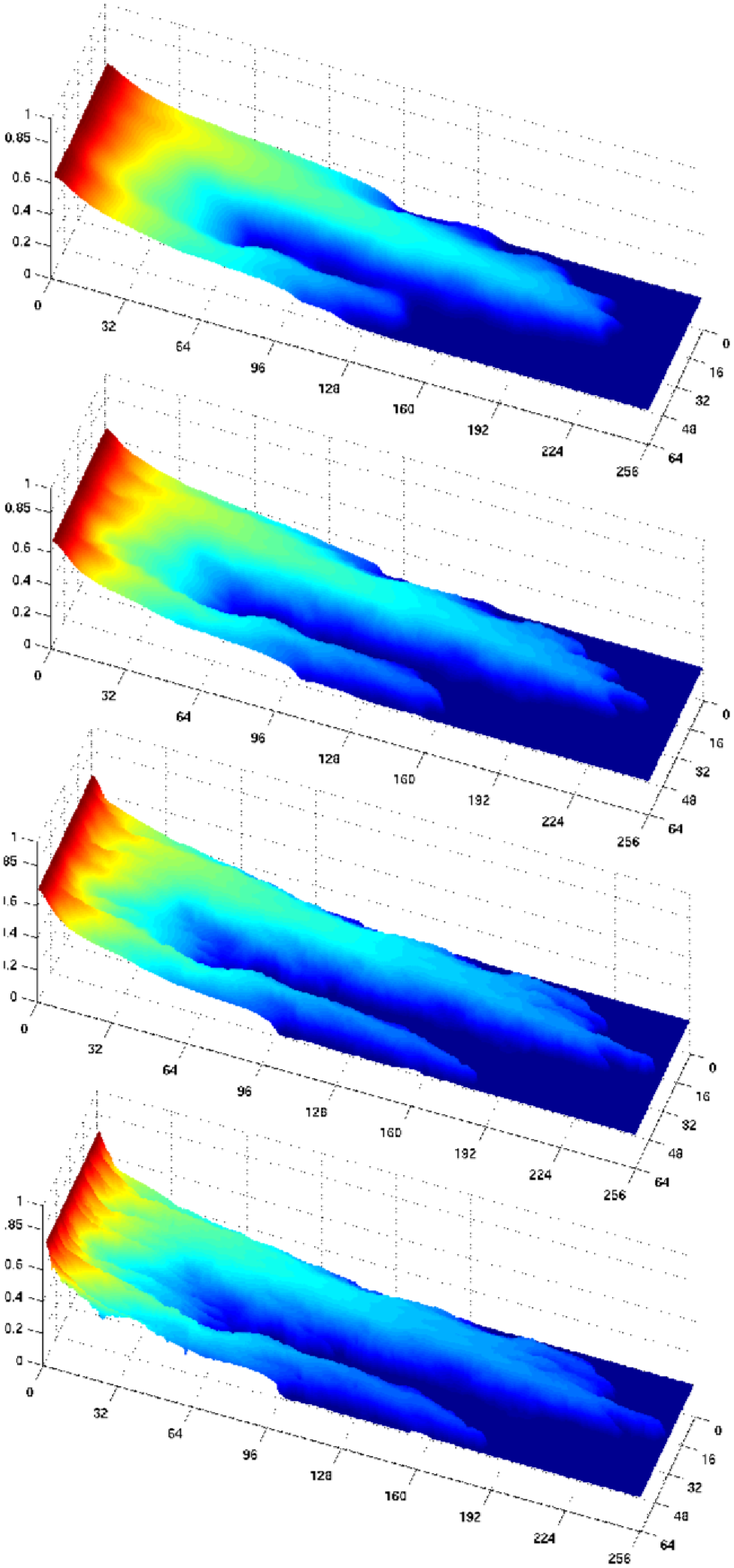}
\caption{Water saturation surface plots
  after $225$ days of simulation in a
  heterogeneous reservoir having $256$ m $\times$ $64$ m,
  with $CV=2.2$ and viscosity
  ratio 20. The first three pictures from top to bottom used
  the NT scheme with grids
  having $256 \times 64$, $512 \times 128$ and $1024 \times
  256$ cells, repectively. The bottom picture shows the KT scheme with a
  grid of $256 \times 64$ cells.}
\label{slab_CV2}
\end{figure}

Figures
\ref{slab}, \ref{slab_CV1}, and \ref{slab_CV2} refer to a
comparative study of the NT and KT schemes, showing the
water saturation surface plots after $275$, $250$ and $225$
days of simulation  for the three different CV values ($CV=0.5, \, ,1.2, \,$ and $2.2$). The results obtained with the NT scheme were computed 
using three computational grid: the coarsest grid with $256 \times 64$ cells, 
and two levels of refinement denoted by NTr and NTrr with $512\times 128$ and $1024 \times 256$ cells, respectively (See the first three pictures from top to bottom of
Figures \ref{slab}, \ref{slab_CV1}, and \ref{slab_CV2}). At the same time,
the bottom pictures in those figures are the results presented by the KT 
dimension by dimension scheme on the coarsest computational grid
of $256 \times 64$ cells. 

\begin{figure}[htbp!]
\includegraphics[scale=0.35]{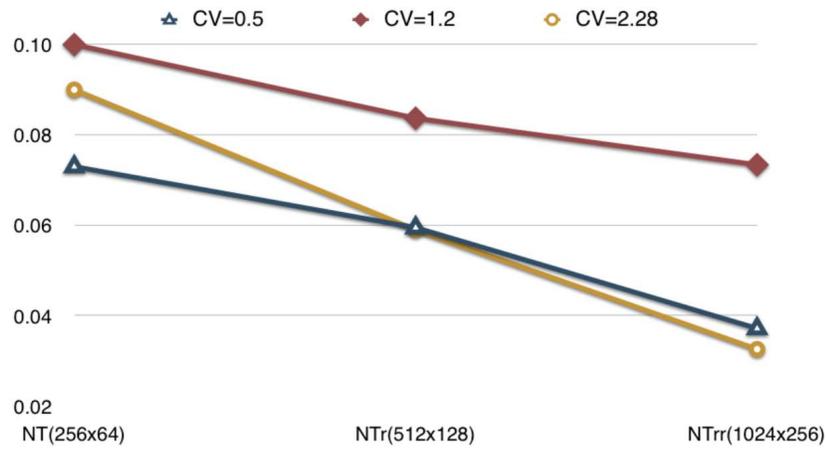}
\caption{The numerical differences in L$\!^2$ norm between the solution of the KT scheme and
the solutions of the NT scheme using three computional grids. As we refine the grid of NT scheme, the differences become smaller.}
\label{fig:graph_erros_numericos}
\end{figure}

For each heterogeneity, we computed the difference 
between the results produced by the NT scheme in the three computational grids
with the corresponding result produced by the KT scheme in the coarsest grid.  We consider the solution
of the NT scheme in the finer grid as the reference solution. The differences are
then computed using the L$\!^2$ norm relative to this reference solution as follows
\begin{equation}
\mbox{numerical differences} = \frac{\parallel{\mbox{F} - \mbox{G}}\parallel_2}{\parallel{\mbox{NTrr}}\parallel_2}.
\end{equation}
Here $F$ stands for the KT solution and $G$ for a NT solution. 
The graph in Figure \ref{fig:graph_erros_numericos} shows these 
differences. Note that as we refine the solution of the NT scheme, 
the differences become smaller indicating a comparable accuracy for the simulations produced by the KT scheme in the coarsest grid and the result produced by the NT scheme 
in the finest grid. These differences indicate  
that one has to refine twice the
grid used in the NT scheme   to produce an 
equivalent solution to the one produced by the KT scheme using the coarsest
grid.

\begin{figure}[htbp]
\centering
\subfigure[\small NT: $64 \! \times \! 64$ grid
]{\label{cfl0_75a}
\includegraphics[scale=0.4]{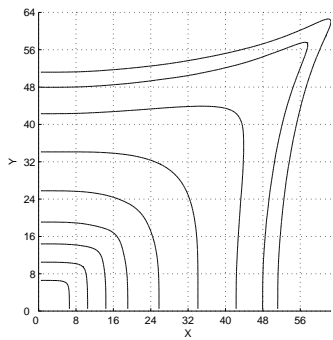}} \hspace{-1.5cm}
\subfigure[\small KT: $64 \! \times \! 64$ grid
]{\label{cfl0_05a}
\includegraphics[scale=0.35]{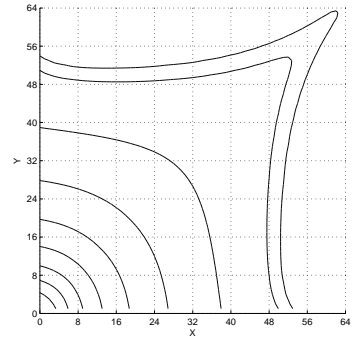}} \\
\vspace{1.5cm}
\subfigure[\small NT: $128 \!\! \times \!\! 128$ grid
]{\label{cfl0_75b}
\includegraphics[scale=0.4]{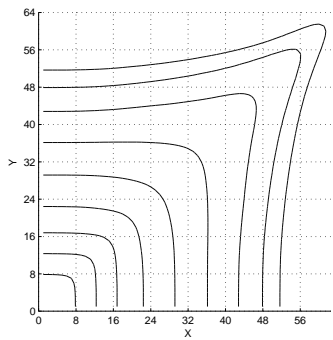}} \hspace{-1.5cm}
\subfigure[\small KT: $128 \! \times \! 128$ grid
]{\label{cfl0_05b}
\includegraphics[scale=0.35]{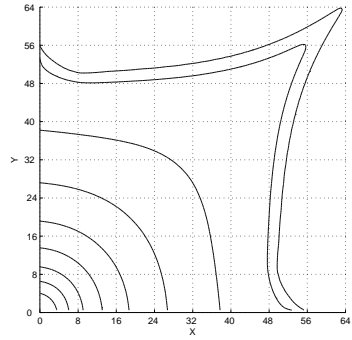}} \\
\caption{Water saturation level curves for two-phase flows in a five-spot well 
configuration - diagonal grid. 
}
\label{5_spot}
\end{figure}

\begin{figure}[htbp]
\centering
\subfigure[\small NT: $90 \! \times \! 90$ grid
]{\label{cfl0_75c}
\includegraphics[scale=0.4]{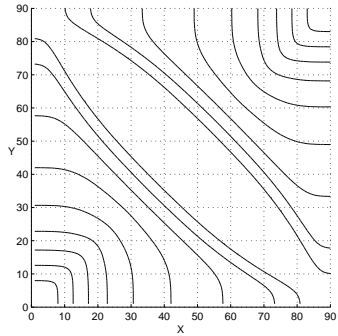}} 
\subfigure[\small KT: $90 \! \times \! 90$ grid
]{\label{cfl0_05c}
\includegraphics[scale=0.35]{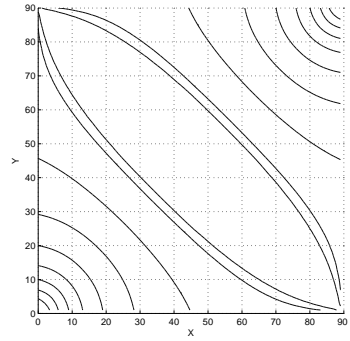}} \\
\vspace{1.5cm}
\subfigure[\small NT: $180 \!\! \times \!\! 180$ grid
]{\label{cfl0_75d}
\includegraphics[scale=0.4]{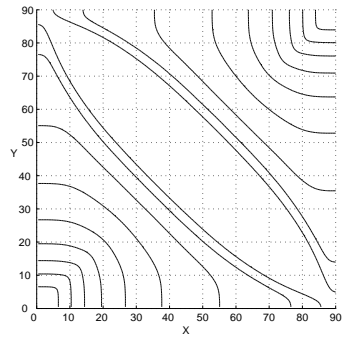}}
\subfigure[\small KT: $180 \! \times \! 180$ grid
]{\label{cfl0_05d}
\includegraphics[scale=0.35]{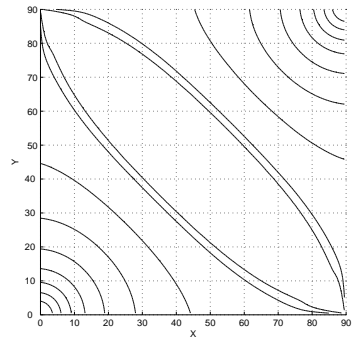}} \\
\caption{Water saturation level curves for two-phase flows in a five-spot well
configuration - parallel grid. 
}
\label{parallel_grid}
\end{figure}

We now turn to the discussion of the set of simulations
performed in a five-spot pattern. In case of a five-spot flood
discretized by a diagonal grid  (Figure \ref{5_spot}), injection takes place
at one corner and production at the diametrically opposite corner;
no flow is allowed across the entirety of the boundary. In case of a five-spot flood discretized by a parallel grid  (Figure \ref{parallel_grid}), injection takes place at two opposite corners (say, bottom left and top right), and production is
through the remaining two corners (say, bottom right and top left). Figures \ref{5_spot} (diagonal grid) and \ref{parallel_grid} (parallel grid) show the saturation level curves after $260$ days of simulation obtained with the NT and KT schemes for two levels of spatial discretization. 

In both Figures \ref{5_spot} and
\ref{parallel_grid}, the pictures on the
left column are the results obtained with the NT scheme and
the ones on the right were computed with the KT scheme. In
these Figures, the grids are refined from top to bottom and
have $64 \, \times \, 64$ and $128 \, \times \, 128$ cells
in the diagonal pattern and $90 \, \times \,
90$ and $180 \, \times \, 180$ cells in the parallel grid.

It is clear that the KT scheme (right column pictures in
Figure~\ref{5_spot} and in Figue~\ref{parallel_grid}) is producing incorrect
boundary behavior. Moreover, as the computational grid is
refined (right column and bottom picture in Figures
\ref{5_spot} and \ref{parallel_grid}) this problem seems to be emphasized.

\clearpage
\clearpage

\pagebreak

\clearpage

\section{Conclusions}\label{sec4}

In one spatial dimension
the KT scheme  is a small modification of the NT scheme which
uses more precise information about the local speed of propagation. 
This approach leads to a very simple numerical recipe producing numerical solutions more accurate then those provided by the NT scheme. On the one hand, In two spatial dimensions the KT scheme uses  numerical fluxes in the $x$ and $y$ directions that can be viewed as generalizations of
one-dimensional numerical fluxes. This is called the dimension by dimension approach.
The NT, in the other hand, uses a genuinely two-dimensional configuration. In the case of a slab geometry, the fluid flows mostly in one direction. For this reason, this flow may be viewed as a one-dimensional flow and the KT scheme is expected to produce very accurate solutions like those presented in Figures \ref{slab}, \ref{slab_CV1} and \ref{slab_CV2}.
In the five-spot problem, the fluid flows in both $x-$ and $y-$directions, causing a genuinely two-dimensional displacement. 
The KT scheme produces
incorrect boundary  behaviors in the
five-spot numerical examples. We remark that this
incorrect behavior is not present in the results produced by the NT
scheme. 
The dimension by dimension approach of the KT scheme might be a source of 
numerical errors for a class of problems with an intrinsic two-dimensional geometry.  These numerical errors may lead to incorrect behavior like those in the five-spot problem.
The authors are currently working on
an improvement of these schemes in order to compute more precisely a genuinely two-dimensional numerical flux. 

\vspace{+1cm}

\paragraph{Acknowledgements.}
The authors wish to thank F. Furtado (University of Wyoming, USA) for several enlightening
discussions and many suggestions during the preparation of this work. S. Ribeiro wishes to thank CAPES/Brazil for a Doctoral fellowship and also F. Pereira and the University of Wyoming for  supporting a Pos-Doctoral study.

\medskip
Received September 8, 2006; accepted February 9, 2007.
\end{document}